\documentclass[11pt,reqno]{amsart}
\usepackage{amssymb}
\usepackage{amscd}
\usepackage[latin1]{inputenc}
\usepackage{amsmath}
\usepackage{mathtools}
\usepackage{tikz}
\usepackage{amsfonts}
\usepackage{faktor}
\usepackage{latexsym}
\usepackage{verbatim}
\usepackage{graphicx}
\usepackage{epsfig}

\usepackage{enumitem}

\usepackage[linktocpage=true]{hyperref}
\usepackage[capitalize]{cleveref}
\hypersetup{citecolor = black,colorlinks,linkcolor = black,urlcolor = black}

\newtheorem{theorem}{Theorem}[section]

\newtheorem{lemma}[theorem]{Lemma}

\newtheorem{corollary}[theorem]{Corollary}

\usepackage[hmargin=3cm,vmargin=3cm]{geometry}

\theoremstyle{definition}
\newtheorem{definition}[theorem]{Definition}
\newtheorem{example}[theorem]{Example}

\theoremstyle{plain}
\newcounter{MainTheoremCounter}

\begin{document}

\title{Spectral gaps and mixing rates of primitive matrices}

\author{Solly Coles}
\address{Solly Coles\\
Department of Mathematics\\
Tufts University\\
177 College Ave.\\
Medford, MA 02155}
\email{solly.coles@tufts.edu}
\urladdr{https://sites.google.com/view/sollycoles}

\author{Ronnie Pavlov}
\address{Ronnie Pavlov\\
Department of Mathematics\\
University of Denver\\
2390 S. York St.\\
Denver, CO 80208}
\email{ronnie.pavlov@du.edu}
\urladdr{http://www.math.du.edu/$\sim$rpavlov/}

\thanks{The second author gratefully acknowledges the support of a Simons Foundation Collaboration Grant.}

\renewcommand{\subjclassname}{MSC 2020}

\begin{abstract}
We obtain lower bounds on spectral gap/upper bounds for mixing rate for arbitrary $n$-state (not necessarily reversible) primitive Markov chains which have exponent $-O(n)$ rather than $-O(n^2)$, a significant improvement on the order obtained from classical arguments. Our improvements come from using the irreducibility distance and chainable index (which are at most linear in $n$) rather than the primitivity index (which can be near $n^2$). This also yields corresponding bounds for arbitrary nonnegative primitive matrices.
\end{abstract}

\maketitle

\section{Introduction}
 
This work concerns mixing rates for primitive matrices and associated Markov chains.
Recall that an $n\times n$ nonnegative matrix $A$ is \textbf{primitive} if there exists a positive integer $R$ (the minimum such $R$ is called the \textbf{primitivity index} of $A$) such that all entries of $A^R$ are positive. By the Perron-Frobenius theorem, any such $A$ has a unique maximal (Perron) eigenvalue $\lambda>0,$ with strictly positive left and right eigenvectors $\ell, r$ (which we assume are normalized so that $\ell r=1$), and the powers of $A/\lambda$ converge to $r \ell$. This motivates the following definition.

\begin{definition}
    The \textbf{(asymptotic) mixing rate} of a primitive matrix $A$ is
    \[
    \eta\coloneq\lim_{m \rightarrow \infty} \|(A/\lambda)^m - (r \ell) \|^{1/m}.
    \]
    
    The \textbf{(normalized additive) spectral gap} of a primitive matrix $A$ is 
    \[
    g = 1 - \frac{|\rho|}{\lambda},
    \]
    where $\rho$ is the eigenvalue of second largest modulus.
\end{definition}

In fact (see for instance Theorem 3.1 of \cite{IL}), it is always the case that $g = 1 - \eta$.
We consider the following fairly elementary question: given an $n \times n$ primitive matrix $A$, what general upper bounds can be given on its mixing rate? Equivalently, what general lower bounds can be given on its spectral gap? 


One common way of obtaining such bounds is by first converting to a primitive Markov chain. Namely, the values $\eta$ and $g$ are unchanged if we replace $A$ by the (Parry) Markov chain $P$ associated to $A$, defined by:
    \[
    p_{ij} = \frac{r_j a_{ij}}{\lambda r_i}.
    \]
Then $P$ is a Markov chain with Perron eigenvalue $1$ and stationary distribution $\pi$ given by $\pi_i= \ell_i r_i$. Also,  $\lambda P=D^{-1}AD$, where $D$ is diagonal and filled with the right Perron eigenvector, implying that $A$ and $P$ have the same values of $\eta$ and $g$. Since $P$ has nonzero entries in the same locations as $A$, it is primitive with the same primitivity index $R$. 
\color{black}

The most standard upper bound on mixing rate for nonreversible Markov chains involves a standard induction argument (see Theorem 4.9 in \cite{LPW}), which yields mixing rate
\[
\eta \leq \sqrt[R]{1 - \delta},
\]
where $R$ is the primitivity index, and 
$\delta = \min p^{(R)}_{ij}/\pi(j)$ over states $i,j$. A tight upper bound on $R$ is $(n-1)^2 + 1$, due to Wielandt (see \cite{schneider}), so this argument yields an upper bound on mixing rate on the order of $1 - \delta/n^2$. 

Choose $i,j$ achieving the minimum in the definition of $\delta$. Due to telescoping for powers of $P$, we have
\[
\delta \geq p^{(R)}_{ij} = \frac{r_j a^{(R)}_{ij}}{\lambda^R r_i} \geq (a/\lambda)^R \frac{r_j}{r_i},
\]
where $a$ is the smallest nonzero entry of $A$. Finally, the shortest path from $j$ to $i$ has length $L$ less than $n$, and so by the eigenvector equation, we have
\[
a^L r_i \leq \lambda^L r_j \Longrightarrow \frac{r_j}{r_i} \geq (a/\lambda)^{n-1}.
\]
Putting this together, we see that this argument yields an upper bound on the mixing rate for $A$ on the order of 
\[
1 - n^{-2} (a/\lambda)^{R+n-1} \leq 1 - n^{-2} (a/\lambda)^{n^2}.
\]

However, this bound seems nonoptimal; its use of Wielandt's bound causes a quadratic factor of $n^2$ in the exponents, and even the conjectured worst-case examples 
have spectral gap on the order of $(a/\lambda)^{O(n)}$ (see Proposition 2.4 of \cite{stanic}).  

By using the idea behind the pseudospectral gap of Paulin (\cite{paulin}), which was itself a refinement of the multiplicative reversibilization bounds for nonreversible Markov chains due to Fill (\cite{fill}), along with recent upper bounds on chainable index (\cite{AGS}), we can give significantly improved lower bounds on spectral gap. Our bounds have exponents on the same order as the \emph{double kite} examples, conjectured by Stani\'c (\cite{stanic}) to have minimal additive spectral gap $\lambda_1 - \lambda_2$ among connected simple graphs for a fixed $n$. The conjecture was recently proved for sufficiently large $n$ by Liu, Tait and Wang (\cite{LTW}). 

The key improvement is that our bounds are expressible in terms of the irreducibility distance $k$ and chainable index $t$, both of which are bounded above by $n$, rather than the possibly much larger primitivity index $R$.

\begin{definition}
    The \textbf{irreducibility distance} $k$ is the minimal positive integer for which all entries of  $A + A^2 + \cdots + A^k$ are strictly positive. The \textbf{chainable index} $t$ is the minimal positive integer for which $A^t (A^T)^t$ is irreducible.
\end{definition}

\color{black}
Our main general bound on the mixing rate/spectral gap is as follows.

\begin{theorem}\label{mainthm}
For any $n \times n$ primitive (nonnegative) matrix $A$, normalized to have largest entry $1$ and smallest nonzero entry $a$, the asymptotic mixing rate of $A$ satisfies
\[
\eta < 1 - \frac{1}{4n^2 t} (a/\lambda)^{8k+4t},
\]
where $\lambda > 0$ is the Perron eigenvalue, $k$ is the irreducibility distance, and $t$ is the chainable index. In addition, the spectral gap is
\[
g = 1 - \eta > \frac{1}{4n^2 t} (a/\lambda)^{8k+4t}.
\]
\end{theorem}


In fact, recent work of Alpin, Guterman, and Shafeev (\cite{AGS}) shows that $t < n$, and trivially $k \leq n$. So this already yields the following corollary:

\begin{corollary}\label{maincor}
    For any $A$ as in \Cref{mainthm}, the asymptotic mixing rate $\eta$ satisfies
    \[
    \eta < 1 - \frac{1}{4n^3} (a/\lambda)^{12n}
    \]
    and the spectral gap satisfies
    \[
    g = 1 - \eta >  \frac{1}{4n^3} (a/\lambda)^{12n}.
    \]
\end{corollary}

We also prove that $t \leq 2k$ (Lemma~\ref{t2k}), yielding a bound purely in terms of the diameter $k$, which could prove useful for large matrices with irreducibility distance $k << n$.
\[
g > \frac{1}{8k^3} (a/\lambda)^{18k}.
\]
In the case of primitive Markov chains, we can also give upper bounds on mixing time:
\[
t_{\mathrm{mix}}(1/4) \leq 2n^3 p^{-12n} \ln(4np^{-2n}) = -4 (\ln p) n^4 p^{-12n} (1 + o(1)) \textrm{ and}
\]
\[
t_{\mathrm{mix}}(1/4) \leq 4k^3 p^{-18k} \ln (4kp^{-3k}) = -12 (\ln p) k^4 p^{-18k} (1 + o(1)).
\]

Our results are similar in spirit to those obtained by Aksoy, Chung, Tait and Tobin (\cite{ACTT}) dealing with the spectral gap of the normalized Laplacian, where they proved an asymptotically sharp lower bound of $(54 + o(1))/n^3$ conjectured by Aldous and Fill. In fact, the question of finding optimal bounds on the `adjacency spectral gap' was posed in Section 4 of their paper, along with other related problems.

\section{Preliminaries on Cheeger bounds and reversibilizations}\label{cheeger}

For a reversible primitive Markov chain $P$, one of the most well-known bounds on spectral gap is Cheeger's inequality, which in some sense measures the smallest bottleneck in the associated directed graph.

\begin{definition}
For any reversible primitive Markov chain $P$, define its \textbf{Cheeger constant} by
\[
h = h(P) = \min_{0 < \pi(S) \leq 1/2} \left\{\frac{\sum_{i \in S} \sum_{j \in S^c} \pi(i) p_{ij}}{\sum_{i \in S} \pi(i)}\right\}.
\]
\end{definition}

If we denote by $\lambda_2$ the second largest eigenvalue of $P$, then Cheeger's inequality states
\[
\frac{h^2}{2} \leq 1 - \lambda_2 \leq 2h.
\]

Technically none of this applies when $P$ is nonreversible, but in \cite{fill}, Fill used so-called reversibilizations to give some relevant bounds. Specifically, one can define the time-reversed version of $P$, denoted $P^*$, by
\[
p^*_{ij} = \frac{\ell_j a_{ji}}{\lambda \ell_i}, \quad \pi_i = \ell_i r_i.
\]
and the \textbf{multiplicative reversibilization} of $P$ by $M = M(P) = P P^*$. It's easily checked that $M$ is reversible and also has stationary distribution $\pi$. 

By (2.27) of \cite{fill}, for any starting distribution $v$ and $m \in \mathbb{N}$, we have
\begin{equation}\label{mixingbound0}
4\|v P^m - \pi\|^2 \leq \chi^2(v \| \pi) \left(1 - h(M)^2/2\right)^m,
\end{equation}
and so $\eta \leq \sqrt{1 - h(M)^2/2}$. 

However, it turns out that $M$ need not even be irreducible, in which case $h(M) = 0$ and this bound is useless. The minimal power where $M(P^j)$ is irreducible is nothing other than the chainable index $t$ of $P$, which is likely much smaller than $R$. In particular, it was shown in \cite{AGS} that $t < n$. It was noted by Paulin in \cite{paulin} that it's possible that powers of $P$ might have reversibilizations with larger spectral gap, which he used in his definition of
\textbf{pseudospectral gap}.

\begin{definition}
    For any $P$, the \textbf{pseudospectral gap} of $P$ is defined as
\[
  g_{ps} = \max_{j \geq 1} \frac{g(P^j (P^*)^j)}{j},
\]
where $g(A)$ is the normalized additive spectral gap of $A$, which can be $0$ if $A$ is not irreducible.
\end{definition}

We do not directly use the pseudospectral gap in this work, but take inspiration from the idea of using powers of $P$ before the multiplicative reversibilization.

\section{Proof of Theorem 1.3}

\begin{proof}
    Suppose that $A,\lambda, k, t$ are as defined in the statement. 
Recall the notation $\ell, r$ for the Perron eigenvectors for $\lambda$ (normalized so that $\ell r = 1$), and the associated Markov Chain $P$ with stationary distribution $\pi$: 
    \[
    p_{ij} = \frac{r_j a_{ij}}{\lambda r_i}, \quad \pi_i = \ell_i r_i.
    \]
\color{black}


For the chainable index $t$, define the multiplicative reversibilization $M$ of 
$P^t$ by $M = P^t (P^*)^t$. 
Then by (\ref{mixingbound0}), for any starting distribution $v$ and $m \in \mathbb{N}$, we have
\begin{equation}\label{mixingbound}
4\|v (P^t)^m - \pi\|^2 \leq \chi^2(v \| \pi) \left(1 - h(M)^2/2\right)^m 
\end{equation}
where $h(M)$ is the Cheeger constant associated to $M$. We now bound $h := h(M)$ from below.
Recall that
\[
h(M) = \min_{0 < \pi(S) \leq 1/2} \left\{\frac{\sum_{i \in S} \sum_{j \in S^c} \pi(i) m_{ij}}{\sum_{i \in S} \pi(i)}\right\}.
\]
Since $M$ is irreducible by definition of $t$ (it has nonzero entries at the same locations as $A^t (A^T)^t$), $h(M) > 0$. In fact $h(M)$ is at least the minimal nonzero entry of $M$ times the minimum of $\pi(i)$ over states $i$.

By definition of $P$, all nonzero entries of $P^t$ contain at least one positive summand of the form
$\frac{r_{i_t} \prod_{q=0}^{t-1} a_{i_q i_{q+1}}}{\lambda^t r_{i_0}}
\geq (a/\lambda)^t \frac{r_{i_t}}{r_{i_0}}$.

For any states $i$ and $j$, there exists a path from $i$ to $j$ of length $d \leq k$ by definition of $k$. By the fact that $r$ is an eigenvector for $A^d$ with eigenvalue $\lambda^d$, we get
\[
\lambda^{d} r_i \geq a^d r_j \Longrightarrow \frac{r_i}{r_j} \geq (a/\lambda)^d \geq (a/\lambda)^k.
\]
Putting this together yields that all nonzero entries of $P^t$ are at least
$(a/\lambda)^{k+t}$, and the same is true for $(P^*)^t$ for similar reasons involving $\ell$.
Finally, for any $i,j$, $\pi(i)/\pi(j) = (\ell_i/\ell_j)(r_i/r_j) \geq (a/\lambda)^{2k}$, implying that $\min \pi(i) \geq (a/\lambda)^{2k}/n$. Therefore, 
\[
h(M) \geq n^{-1} (a/\lambda)^{4k+2t},
\]
and so by (\ref{mixingbound}) we have
\begin{equation}\label{mainbound}
4\|v P^{tm} - \pi\|^2 \leq
\chi^2(v \| \pi) \left(1 - \frac{1}{2n^2} (a/\lambda)^{8k+4t}\right)^m
\end{equation}
for all $m$, implying that the mixing rate of $P$ (which is equal to that of $A$) is bounded from above by
\[
\left(1 - \frac{1}{2n^2} (a/\lambda)^{8k+4t}\right)^{1/2t} \leq 
1 - \frac{1}{4n^2 t} (a/\lambda)^{8k+4t}.
\]
If we denote the second largest eigenvalue modulus of $P$ by $\rho$, then it's well-known that asymptotic mixing rate is just the modulus of the second largest eigenvalue, and so we also get that the spectral gap of $P$ (which is the same as that of $A$) is bounded from below 
by $\frac{1}{4n^2 t} (a/\lambda)^{8k+4t}$.

\end{proof}

As noted earlier, since $k \leq n$ and $t < n$ by \cite{AGS}, we get the following bounds purely in terms of $n$:
\[
    \eta < 1 - \frac{1}{4n^3} (a/\lambda)^{12n}
    \]
    and the spectral gap satisfies
    \[
    g = 1 - \eta >  \frac{1}{4n^3} (a/\lambda)^{12n}.
    \]

\section{Bounds in terms of $k$}
It is desirable to have a bound on $\eta$ purely in terms of $k$ (and $a$ and $\lambda$), since for many natural systems there are good upper bounds on $k$ much less than $n$. 

\begin{lemma}\label{t2k}
For every primitive $A$, $t \leq 2k$.
\end{lemma}

\begin{proof}
Theorem 3.5 from \cite{AGS} shows that the number of irreducible components of $A^m (A^T)^m$ is strictly decreasing in $m$ (until it becomes $1$). Since there were clearly at most $n$ irreducible components for $m = 0$, this yields their result that $t < n$. 

For this lemma, it therefore suffices to show that $A^k (A^T)^k$ has less than or equal to $k+1$ irreducible components, since then in $k$ more steps we must have $A^{2k} (A^T)^{2k}$ with a single component.

However, this is fairly simple. Fix any state $i$, and classify states by minimal distance (with respect to edges from $A$) to $i$; $V_d$ is the set of all states $j$ for which there is a path from $j$ to $i$ of length $d$ (equivalently, 
$(A^d)_{ji} > 0$), but none of length less than $d$. Then, by definition, all states of $V_d$ are connected in $A^k (A^T)^k$; for any pair $a,b \in V_d$, there are paths
of length $d$ from $a$ to $i$ and $b$ to $i$, and then appending an arbitrary path of length $k-d$ starting at $i$ yields paths of length $k$ ending at the same vertex.

Since $V_0, \ldots, V_k$ partition states by definition of $k$, we had at most $k+1$ components in $A^k (A^T)^k$, completing the proof.

\end{proof}

We now have the following corollary of Theorem~\ref{mainthm}.

\begin{corollary}\label{irredcor}
    For any $A$ as in \Cref{mainthm}, the asymptotic mixing rate $\eta$ satisfies
    \[
    \eta < 1 - \frac{1}{8k^3} (a/\lambda)^{18k}.
    \]
    and the spectral gap satisfies
    \[
    g = 1 - \eta > \frac{1}{8k^3} (a/\lambda)^{18k}.
    \]
\end{corollary}

\begin{proof}
We note that since $A/\lambda + \cdots + (A/\lambda)^k$ is a positive matrix, and all entries are at least $(a/\lambda)^{k}$, its Perron eigenvalue is at least $n(a/\lambda)^k$. On the other hand, its Perron eigenvalue by definition is $k$. This implies
\begin{equation}\label{ntok}
n \leq k (a/\lambda)^{-k},
\end{equation}
and substituting this and $t \leq 2k$ into \cref{mainthm} completes the proof.
\end{proof}

\section{Mixing time}
Finally, specifically in the case of primitive Markov chains, we can obtain associated upper bounds on mixing time.

\begin{theorem}
For any primitive Markov chain on $n$ states with minimum nonzero entry $p$ and irreducibility distance $k$, the mixing time 
$t_{\mathrm{mix}}(1/4)$ satisfies
\[
t_{\mathrm{mix}}(1/4) \leq 2n^3 p^{-12n} \ln(4np^{-2n}) = -4 (\ln p) n^4 p^{-12n} (1 + o(1)) \textrm{ and}
\]
\[
t_{\mathrm{mix}}(1/4) \leq 4k^3 p^{-18k} \ln (4kp^{-3k}) = -12 (\ln p) k^4 p^{-18k} (1 + o(1)).
\]
\end{theorem}

\begin{proof}
Assume that $P$ is a primitive Markov chain and recall (\ref{mainbound}):
\begin{equation}
4\|v P^{tm} - \pi\|^2 \leq
\chi^2(v \| \pi) \left(1 - \frac{1}{2n^2} (a/\lambda)^{8k+4t}\right)^m
\end{equation}
For a Markov chain, after normalizing to make largest entry $1$, $a/\lambda$ is simply the smallest nonzero entry of $P$, call it $p$.
Also, clearly $\chi^2(v \| \pi) \leq (\min \pi(i))^{-1}$, and the argument used previously establishes $(\min \pi(i))^{-1} \leq np^{-2k}$. Putting this together yields that if 
\begin{equation}
np^{-2k} \left(1 - \frac{1}{2n^2} p^{8k+4t}\right)^m < 1/4,
\end{equation}
then $tm$ is an upper bound on mixing time. Therefore,
\[
t_{\mathrm{mix}}(1/4) \leq 2n^2 t p^{-(8k+4t)} \ln (4np^{-2k}).
\]
Just using the inequalities $t,k \leq n$ yields
\[
t_{\mathrm{mix}}(1/4) \leq 2n^3 p^{-12n} \ln(4np^{-2n}) = -4 (\ln p) n^4 p^{-12n} (1 + o(1)),
\]
and using $t \leq 2k$ and (\ref{ntok}) yields
\[
t_{\mathrm{mix}}(1/4) \leq 4k^3 p^{-18k} \ln (4kp^{-3k}) = -12 (\ln p) k^4 p^{-18k} (1 + o(1)).
\]
\end{proof}

\section{Examples}
The following class of simple examples demonstrates that spectral gaps of $0$-$1$ matrices can indeed be exponentially small in $n$ and $k$. They are essentially just a directed looped version of the double kite model treated in \cite{stanic}.

\begin{example}
For any $2 \leq m \leq L$, define the digraph consisting of 
two disjoint copies of complete digraphs (with self-loops) on $m$ vertices with distinguished vertices $a\neq b$ in the first and $c\neq d$ in the second, together with vertex-disjoint directed paths of length $L$ from $a$ to $c$ and from $d$ to $b$ which share only $a,b,c,d$ with the complete graphs (i.e. each path contains $L-1$ new vertices); see Figure~\ref{pic}).
The number of vertices is $n = 2m + 2L - 2 \leq 4L$.

\begin{figure}[h]
\centering
\begin{tikzpicture} 
\usetikzlibrary{arrows.meta} 
\draw (3,0) circle (1cm);
\draw (-3,0) circle (1cm);
\node at (-3,0){$K_m$};
\node at (3,0){$K_m$};

\draw[->,thick]  (-2.5,0.5) to[bend left] (2.5,0.5);
\draw[->,thick]  (2.5,-0.5) to[bend left] (-2.5,-0.5);
\draw[fill=black] (-2.5,0.5) circle (0.05cm);
\draw[fill=black] (2.5,0.5) circle (0.05cm);
\draw[fill=black] (-2.5,-0.5) circle (0.05cm);
\draw[fill=black] (2.5,-0.5) circle (0.05cm);
\node at (-2.7,0.5){$a$};
\node at (2.7,0.5){$c$};
\node at (-2.7,-0.5){$b$};
\node at (2.7,-0.5){$d$};

\node at (0,1.5){length $L$};
\node at (0,-1.5){length $L$};

\end{tikzpicture}
\caption{}\label{pic}
\end{figure}

It's not hard to verify that the associated $0$-$1$ matrix $A$ has characteristic polynomial
\[
p(x) = x^{2m-4} \left(x^{2L}(x-m)^2 - 1\right).
\]
Therefore, the Perron eigenvalue $\lambda$ and some other eigenvalue $\rho$ satisfy
\[
\lambda^{L} (\lambda-m) = 1, \quad \rho^{L} (\rho - m) = -1.
\]
Then $\lambda - m = \frac{1}{\lambda^L} < m^{-L}$. In addition, plugging in $m - 2m^{-L}$ to $-x^L(x-m)$ yields
\[
(m - 2m^{-L})^L \cdot 2m^{-L} = 2 \left(1 - \frac{2}{m^{L+1}}\right)^L > 2\left(1 - \frac{2L}{m^{L+1}}\right),
\]
which is at least $1$ since $m^{L+1} \geq 4L$. Therefore, $\rho > m - 2m^{-L}$ and
\[
g \leq 1 - \rho/\lambda < 1 - \frac{m - 2m^{-L}}{\lambda} = (\lambda - m + 2m^{-L})/\lambda < 3m^{-L}/\lambda.
\]
Finally,
\[
(\lambda/m)^{L} < \left(1 + m^{-L-1}\right)^{L} < e^{1/4}
\]
since $m^{L+1} \geq 4L$, and so
\[
g < 3m^{-L}/\lambda < 3e^{1/4} \lambda^{-L-1} < 4 \lambda^{-L-1}.
\]

This gives a spectral gap of less than
$4 \lambda^{-L-1} < 4 (a/\lambda)^{n/4}$. Finally, the irreducibility distance is
$k = 2L + 2$, and so the spectral gap is also less than $4 (a/\lambda)^{k/2}$.

\end{example}

\section{AI declaration}

The authors acknowledge the use of Astra GPT-6 in the preparation of this paper, for proofreading and literature searches, including identifying the relevance of the article \cite{AGS} and its bound on the chainable index. All proofs and writing in the paper were written by the authors, who take full responsibility for its mathematical correctness and final text.
\bibliographystyle{plain}
\bibliography{markov}

\end{document}